\smallskip\documentclass[12pt]{article}
\usepackage{amssymb,amsmath,times}

 \makeatletter
\renewcommand{\@listI}{\setlength{\topsep}{0pt}\setlength{\itemsep}{0pt}}

\renewcommand{\@begintheorem}[2]{\begin{trivlist}\it%
\item[\hspace{\parindent}\hspace{\labelsep}{\bf #1\ #2.}]}
\renewcommand{\@endtheorem}{\end{trivlist}}

\newcounter{theoremm}
\newenvironment{theoremm}{\refstepcounter{theoremm}\begin{trivlist}\it%
\item[\hspace{\parindent}\hspace{\labelsep}%
{\bf Theorem\ \thetheoremm{\boldmath $'$}.}]}{\end{trivlist}}

\makeatother

\makeatletter
\renewcommand{\@biblabel}[1]{\hfill #1.}
\makeatother

\newcommand{\Var}{\mathop{\rm Var}}
\newcommand{\Cov}{\mathop{\rm Cov}}

\newcommand{\myp}{\mbox{$\:\!$}}
\newcommand{\mypp}{\mbox{$\;\!$}}

\def\RR{\mathbb R}

\def\ZZ{\mathbb Z}

\newtheorem{theorem}{Theorem}
\newtheorem{lemma}{Lemma}
\newtheorem{remark}{Remark}

\begin{document}
%{\large\centerline{\underline{MATHEMATICS}}}
%
%\smallskip
%\noindent
%UDC 519.21
%\null
\medskip
\smallskip
\begin{center}
{\large\bf CENTRAL LIMIT THEOREM FOR RANDOM PARTITIONS\\[.4pc] UNDER THE
PLANCHEREL MEASURE}
\\[18pt]
{\large\bf L.\,V.\ Bogachev${}^*$ and Z.~G.\ Su${}^{**}$ }\\[6pt]
%\large
$^*$\,{Department of Statistics, University of Leeds,
United
Kingdom.\\
E-mail: bogachev@maths.leeds.ac.uk}\\[4pt]
${}^{**}$\,{Department of Mathematics, Zhejiang University,
Hangzhou, P.R.~China.\\ E-mail: suzhonggen@zju.edu.cn}
\end{center}

%
% 60F05 Central limit and other weak theorems
%
% 60C05 Combinatorial probability
% 05A17 Partitions of integers
% 05E10 Tableaux, representations of the symmetric group

%\large

\bigskip
\hfill{\it To the memory of Sergei Kerov (1946--2000)}

% Abstract as submitted to arXiv.math
%
%In this work, we obtain the central limit theorem for fluctuations of Young
%diagrams around their limit shape in the bulk of the "spectrum" of
%partitions of a large integer n (under the Plancherel measure).
%More specifically, we show that, under the suitable normalization
%(growing like the square root of log n, the corresponding random
%process converges, in the sense of finite dimensional
%distributions, to a Gaussian process with independent values. The
%proof uses heavily the determinantal structure of the correlation
%functions and is based on the application of the
%Costin--Lebowitz--Soshnikov central limit theorem. At the spectrum
%edges, the fluctuation asymptotics is expressed in terms of the
%largest members of the Airy ensemble; in particular, at the upper
%edge the limit distribution appears to be discrete (without any
%normalization). These results admit an elegant symmetric
%reformulation under the rotation of Young diagrams by 45 degrees,
%where the normalization no longer depends on the location of the
%spectrum point. We also discuss the link of our central limit
%theorem with an earlier result by S.V. Kerov on the convergence to
%a generalized Gaussian process.

\bigskip
\begin{center}
\large 1. INTRODUCTION
\end{center}

\smallskip
A \emph{partition} of a natural number $n$ is any integer sequence
$\lambda=(\lambda_1,\lambda_2,\dots)$ such that
$\lambda_1\ge\lambda_2\ge\dots\ge0$ and
$\lambda_1+\lambda_2+\dots=n$ (notation: $\lambda\vdash n$). In
particular, $\lambda_1=\max\{\lambda_i\in\lambda\}$. Every
partition $\lambda\vdash n$ can be represented geometrically by a
planar shape called the \emph{Young diagram}, consisting of $n$
unit cell arranged in consecutive columns, containing
$\lambda_1,\lambda_2,\dots$ cells, respectively.

On the set ${\mathcal P}_n:=\{\lambda\vdash n\}$ of all partitions
of a given $n$, consider the \emph{Plancherel measure}
\begin{equation}\label{eq:Plancherel}
P_n(\lambda):=\frac{d^2_\lambda}{n!}\,, \qquad \lambda
\in{\mathcal P}_n,
\end{equation}
where $d_\lambda $ is the number of standard tableaux of a given
shape $\lambda$, that is, the total number of arrangements of the
numbers $1,\dots,n$ in the cells of the Young diagram $\lambda\in
{\mathcal P}_n$, such that the numbers increase in each row (from
left to right) and each column (bottom up). Note that $d_\lambda$
also equals the dimension of the irreducible (complex)
representation of the symmetric group ${\mathfrak S}_n$ (i.e., the
group of permutations of order $n$), indexed by the partition
$\lambda$ (see \cite{F,VK,VK1}). According to the
RSK\footnote{Robinson--Schensted--Knuth.} correspondence (see
\cite{F}), any permutation $\sigma\in{\mathfrak S}_n$ is
associated with exactly one (ordered) pair of standard tableaux of
the same shape $\lambda\in{\mathcal P}_n$. Since there are $n!$
such permutations, this implies the Burnside identity (see
\cite{VK1})
$$
\sum_{\lambda\in{\mathcal P}_n} d^2_\lambda=n!,
$$
thus the measure $P_n$ defined in (\ref{eq:Plancherel}) determines
a probability distribution on ${\mathcal P}_n$.

The Plancherel measure arises naturally in
representation-theoretic, combinatorial, and probabilistic
problems (see \cite{D}). For example, the RSK correspondence
implies that the largest term $\lambda_1$ of the partition
$\lambda$ associated with a given permutation $\sigma\in{\mathfrak
S}_n$ equals the length $\ell_n$ of the longest increasing
subsequence contained in $\sigma$. Therefore, the Plancherel
distribution of $\lambda_1$ coincides with the distribution of
$\ell_n$ in a random (uniformly distributed) permutation
$\sigma\in S_n$ (see \cite{BDJ}).\footnote{The interest in the
asymptotic behavior of the random variable $\ell_n$ was stimulated
by the Ulam problem (see \cite{D}). The problem was settled by
Vershik and Kerov \cite{VK} who showed that $\ell_n/\sqrt{n}\to 2$
in probability (cf.\ formula (\ref{eq:lambda1}) below).}

The upper boundary of the Young diagram corresponding to the
partition $\lambda\in{\mathcal P}_n$ can be viewed as the graph of
a stepwise (left-continuous) function $\lambda(x)$, $x\ge0$,
defined by
\begin{equation}\label{eq:lambda}
\lambda(x):=\lambda_0\myp{\bf 1}_{\{0\}}(x)+\sum_{i=1}^\infty
\lambda_i\myp {\bf 1}_{(i-1,i]}(x)\equiv \lambda_{\lceil x\rceil},
\end{equation}
where ${\bf 1}_B(x)$ is the characteristic function (indicator) of
set $B$ and $\lceil x\rceil:=\min\{m\in\ZZ: m\ge x\}$ is the
ceiling integer part of $x$. Logan and Shepp \cite{LS} and,
independently, Vershik and Kerov \cite{VK} (see also \cite{VK1})
have discovered that, as $n\to\infty$, a typical Young diagram,
suitably scaled, has a ``limit shape'' represented by the graph of
some function $y=\omega(x)$. This means that for the overwhelming
majority of partitions $\lambda\in{\cal P}_n$ (with respect to the
Plancherel measure $P_n$), the boundary of their scaled Young
diagrams is contained in an arbitrarily small vicinity of the
graph of $\omega(x)$.

More specifically, set
\begin{equation}\label{eq:tilde}
\bar\lambda_n(x):=\frac{1}{\sqrt n}\,\lambda(\sqrt n x),\qquad
x\ge0,
\end{equation}
and consider the function $y=\omega(x)$ defined by the parametric
equations
\begin{equation}\label{eq:omega}
x=\frac{2}{\pi}\,(\sin\theta-\theta\cos\theta),\qquad
y=x+2\cos\theta,\qquad 0\le \theta\le \pi.
\end{equation}
The function $\omega(x)$ is decreasing on $[0,2]$ and
$\omega(0)=2$, \,$\omega(2)=0$. Define $\omega(x)$ as zero for all
$x>2$. Then the random process $\bar\lambda_n(x)$ satisfies the
following law of large numbers \cite{VK,VK1,LS}:
$$
\forall\myp \varepsilon>0,\quad
\lim_{n\to\infty}P_n\Bigl\{\sup_{x\ge 0}
\bigl|\myp\bar{\lambda}_n(x)-\omega(x)\bigr|>\varepsilon\Bigr\}=0.
$$
In particular, for $x=0$ it follows that the maximal term in a
typical partition asymptotically behaves like $2\sqrt{n}$:
\begin{equation}\label{eq:lambda1}
\forall\myp \varepsilon>0,\quad
\lim_{n\to\infty}P_n\Bigl\{\Bigl|\myp\frac{\lambda_1}{\sqrt{n}}
-2\Bigr|> \varepsilon\Bigr\}=0.
\end{equation}

\begin{remark}\label{rm1}
\normalfont Due to the invariance of the Plancherel measure under
the transposition of Young diagrams
$\lambda\leftrightarrow\lambda'$ (when the columns of the diagram
$\lambda$ become rows of the transposed diagram $\lambda'$ and
vice versa), the same law of large numbers holds for
$\lambda'_1=\#\{\lambda_i\in\lambda\}$ (i.e., for the number of
terms in the random partition $\lambda$).
\end{remark}

A natural question about fluctuations of the random function
${\bar\lambda}_n$ around the limit curve $\omega$ was posed in
\cite{LS} (see also \cite{VK1}), but it remained open for more
than 15 years. Kerov \cite{K} (see also \cite{IO}) gave a partial
answer by establishing the convergence of the random process
\begin{equation}\label{eq:Delta}
\Delta_n(x):=\sqrt{n}\,\bigl(\bar\lambda_n(x)-\omega(x)\bigr)
\end{equation}
to a generalized Gaussian process (without any further
normalization!). To state this result more precisely, it is
convenient to pass to the coordinates $u=x-y$, $v=x+y$, which
corresponds to anticlockwise rotation by $45^\circ$ and dilation
by $\sqrt{2}$. In the new coordinates, the boundary of the scaled
Young diagram is determined by the piecewise linear (continuous)
function $\widetilde\lambda_n(u)$ and the limit shape is given by
$$
  \Omega(u):=\left\{\begin{array}{ll}
  \displaystyle \frac{2}{\pi}\left(u\arcsin \frac{u}{2}
  +\sqrt{4-u^{2}}\,\right),&|u|\le 2,\\[.5pc]
  |u|,& |u|\ge 2.
  \end{array}\right.
$$
Then, according to \cite{K,IO}, the random process
\begin{equation}\label{eq:tilde_Delta}
\widetilde\Delta_n(u):=
\sqrt{n}\,\bigl(\,\widetilde\lambda_n(u)-\Omega(u)\bigr)
\end{equation}
converges in distribution to a generalized Gaussian process
$\widetilde\Delta(u)$, $u\in[-2,2]$, defined by the
formal random series
\begin{equation}\label{eq:genDelta}
 \widetilde\Delta(u)|_{u=2\cos\theta}
 =\frac{2}{\pi}\sum_{k=2}^\infty
\frac{X_k \sin(k\theta)}{\sqrt{k}}\,,\qquad \theta\in[0,\pi],
\end{equation}
where $X_2,X_3,\dots$ are independent random variables with
standard normal distribution ${\mathcal N}(0,1)$. Convergence of
$\widetilde\Delta_n(\cdot)$ to $\widetilde\Delta(\cdot)$ is
understood in the sense of generalized functions. It is convenient
to choose test functions in the form of the modified Chebyshev
polynomials of the second kind (see \cite{IO}), defined by the
formula
\begin{equation}\label{eq:U}
  U_{k}(u)|_{u=2\cos\theta}:=
  \frac{\sin((k+1)\theta)}{\sin\theta}\,,\qquad
\theta\in[0,\pi].
\end{equation}
Then one can show (see details in \cite{IO}) that for
$k=2,3,\dots$,
\begin{equation}\label{eq:genlim}
\begin{aligned}
\int_{\RR} \widetilde\Delta_n(u)\,U_{k-1}(u)\,du
&\stackrel{d}{\longrightarrow}
\int_{-2}^2 \widetilde\Delta(u)\,U_{k-1}(u)\,du\\
&=2\int_{0}^\pi
\widetilde\Delta(2\cos\theta)\sin(k\theta)\,d\theta =
\frac{2X_k}{\sqrt{k}}
\end{aligned}
\end{equation}
(here and below, the symbol $\stackrel{d\,}{\to}$ denotes
convergence in distribution).

\begin{remark}\label{rm2}
\normalfont Note that $U_0(u)\equiv 1$ (see (\ref{eq:U})), in
which case we have
$$
\int_{\RR}\widetilde\Delta_n(u)\,U_0(u)\,du
=\int_{\RR}\widetilde\Delta_n(u)\,du= 2\int_{0}^\infty
\Delta_n(x)\,dx=0,
$$
because the area of the scaled Young diagram equals $1$, as well
as the area under the graph of $y=\omega(x)$. This explains why
$k\ge 2$ in (\ref{eq:genDelta}) and (\ref{eq:genlim}). Also note
that fluctuations outside $[-2,2]$ are negligible, since
$\int_{|u|>2}\varphi(u)\,\widetilde\Delta_n(u)\,du\to 0$ in
probability for any test function $\varphi$ with compact support
\strut(see \cite{IO}).
\end{remark}

\begin{remark}\label{rm3}
\normalfont A similar result about convergence to a generalized
Gaussian process for eigenvalues of random matrices in the
Gaussian Unitary Ensemble (GUE) was obtained by Johansson
\cite{J0} (see further discussion of these results in~\cite{IO}).
\end{remark}

However, a ``localized'' version of the central limit theorem for
random partitions (i.e., for fluctuations at a given point) has
not been known as yet. On the one hand, the existence of such a
theorem would have seemed quite natural (at least, in the bulk of
the partition ``spectrum''\footnote{We use the term ``spectrum''
informally by analogy with the GUE, to refer to the variety of
partition's terms $\lambda_i\in\lambda$ (cf.\ the book \cite{ABT}
where this term is used in a general context of combinatorial
structures characterized by their components).}, i.e., for
$\lambda_i\in\lambda\vdash n$ such that $i/n\sim x\in(0,2)$); on
the other hand, Kerov's result on generalized convergence might
cast some doubt on the validity of the usual convergence.

Note that the asymptotic behavior of fluctuations at the upper
edge of the limiting spectrum (corresponding to $x=0$) is
different from Gaussian. As was shown in \cite{BDJ} for
$\lambda_1$ and in \cite{BOO,J,O} for any $\lambda_i$ with fixed
$i=1,2,\dots$,
\begin{equation}\label{eq:Airy}
\lim_{n\to\infty}P_n\biggl\{\frac{\lambda_i-2\sqrt n}{n^{1/6}}\le
z\biggr\} =F_{i}(z),\qquad z\in\RR,
\end{equation}
where $F_i$ is the distribution function of the $i$-th largest
point in the so-called Airy random point process, discovered
earlier in connection with the limit distribution of the largest
eigenvalues for random matrices from the GUE (see~\cite{TW}). In
particular, $F_1(\cdot)$ is known as the Tracy--Widom distribution
function.

From the point of view of Kerov's limit theorem (see
(\ref{eq:genlim})), the extreme values $\lambda_1,\lambda_2,\dots$
might present a danger, since according to formula (\ref{eq:Airy})
the fluctuations of the process $\Delta_n(x)$ near $x=0$ are large
(of order of $n^{1/6}$). As this theorem shows, the edge of the
spectrum in fact does not give any considerable contribution into
the integral fluctuations. Let us stress, however, that the
situation in the bulk of the spectrum remained unclear.

\bigskip
\begin{center}
\large 2. MAIN RESULTS
\end{center}

\smallskip
In our first result, we establish the central limit theorem for
the random variable $\Delta_n(x)$ given by (\ref{eq:Delta}). Set
\begin{equation}\label{eq:Y}
Y_n(x):= \frac{2\myp\theta_x\myp\Delta_n(x)}{\sqrt{\log
n\mathstrut}}\,,
\end{equation}
where $\theta_x=\arccos\frac{\omega(x)-x}{2}$ is the value of the
parameter $\theta$ in equations (\ref{eq:omega}) corresponding to
the coordinates $x$ and $y=\omega(x)$.

\begin{theorem}\label{th1}
For each\/ $0<x<2$, the distribution of the random variable\/
$Y_n(x)$ with respect to the Plancherel measure\/ $P_n$ converges,
as\/ $n\to\infty$, to the standard normal distribution\/
${\mathcal N}(0,1)$.
\end{theorem}

\begin{remark}\label{rm4}
\normalfont One can show that Theorem \ref{th1} also holds for
$Y_n(x_n)$ if $x_n\to x\in(0,2)$ as $n\to\infty$.
\end{remark}

The local structure of correlations of the random process
$\Delta_n(x)$ is described by the following theorem. We write $c_n\asymp
1$ if $c_n n^{\varepsilon}\to\infty$, $c_n n^{-\varepsilon}\to 0$ for any
$\varepsilon>0$, and $a_n\asymp b_n$ if $a_n/b_n\asymp 1$.

\begin{theorem}\label{th2}
Fix\/ $x_0\in(0,2)$ and let\/ $x_1,\dots,x_m\in(0,2)$ be such
that\/ $|x_0-x_i|\asymp n^{-s_i/2}$, where\/ $0\le s_i\le 1$
\,$(i=1,\dots,m)$. For\/ $i=0$, set formally\/ $s_0=1$. Then the
random vector\/ $(Y_n(x_0),\dots,Y_n(x_m))$ converges in
distribution, as\/ $n\to\infty$, to a Gaussian vector\/
$(Z_{s_0},\dots,Z_{s_m})$ with zero mean and covariance matrix\/
$K$ with the elements\/ $K(s_i,s_i)=1$,
$K(s_i,s_j)=\min\{s_i,s_j\}$ \,$(i\ne j)$.
\end{theorem}

Note that by Theorem \ref{th2}, the covariance between $Y_n(x)$
and $Y_n(x')$ asymptotically decays as the distance $|x-x'|$
grows:
$$
 |x-x'|\asymp n^{-s/2}\quad\Rightarrow\quad
 \lim_{n\to\infty}\Cov\bigl(Y_n(x),Y_n(x')\bigr)=s.
$$
In particular, if $|x-x'|\asymp n^{-1/2}$ (i.e., $s=1$), then
$(Y_n(x),Y_n(x'))\stackrel{d}{\to} (Z_1,Z_1)$, while if $x'$ is at
a fixed distance from $x$ (i.e., $s=0$) then $Y_n(x)$, $Y_n(x')$
are asymptotically independent.

\begin{remark}\label{rm5}
\normalfont Results similar to Theorems \ref{th1} and \ref{th2}
were obtained by Gustavsson \cite{G} for eigenvalues in the bulk
of the spectrum of random matrices in the GUE.
\end{remark}

Let us point out that Theorems \ref{th1} and \ref{th2} can be
reformulated in coordinates $u,v$ (see Sect.~1). To this end, one
needs to find the sliding projection (divided by $\sqrt{2}$\,) of
the deviation $\widetilde\Delta_n(u)$ (see~(\ref{eq:tilde_Delta}))
onto the line $u+v=0$ along the tangent of the graph
$\Omega(\cdot)$ at point $u$. Differentiating equations
(\ref{eq:omega}), we get
$$
\frac{dv}{du}=\frac{x'_\theta+y'_\theta}{x'_\theta-y'_\theta}
=\frac{2\myp\theta_x}{\pi}-1\,,
$$
which implies
$$
\widetilde\Delta_n(u)=\frac{2\myp\theta_x}{\pi}\,\Delta_n(x_n)
\left(1+\eta_n\right),
$$
where $x_n\to x$ and $\eta_n\to0$ (in probability). Hence, setting
\begin{equation}\label{eq:tildeY}
\widetilde Y_n(u):=\frac{\pi\myp\widetilde\Delta_n(u)}{\sqrt{\log
n\mathstrut}}\,,\qquad u\in\RR,
\end{equation}
and using Remark \ref{rm4}, we obtain the following elegant
versions of Theorems \ref{th1} and \ref{th2}, where the
normalization constant is the same for all points.

\begin{theoremm}\label{th1'}
For each\/ $-2<u<2$, the distribution of the random variable
$\widetilde Y_n(u)$ with respect to the Plancherel measure $P_n$
converges, as $n\to\infty$, to the normal distribution ${\mathcal
N}(0,1)$.
\end{theoremm}

\begin{theoremm}\label{th2'}
Let $u_0,\dots,u_m\in(-2,2)$, $|u_0-u_i|\asymp n^{-s_i/2}$, $0\le
s_i\le 1$, with the conventions as in Theorem $\ref{th2}$. Then
the random vector\/ $\bigl(\widetilde Y_{n}(u_0),\dots,\widetilde
Y_{n}(u_m)\bigr)$ converges in distribution, as $n\to\infty$, to a
Gaussian vector $(Z_{s_0},\dots,Z_{s_m})$ with zero mean and the
same covariance matrix $K$.
\end{theoremm}

\begin{remark}\label{rm6}
\normalfont The covariance function $K(s,s')$ of Theorems
\ref{th2} and \ref{th2'}$'$ determines a Gaussian process $Z_s$ on
$[0,1]$ (with zero mean), which can be represented as
$$
  Z_s\stackrel{d}{=}
  W_s+\zeta_s\sqrt{1-s}\,, \qquad 0\le s\le 1,
$$
where $W_s$ is a standard Wiener process and $\{\zeta_s,\ 0\le
s\le 1\}$ is a family of mutually independent random variables
(also independent of $W_s$) with normal distribution ${\mathcal
N}(0,1)$. This decomposition shows that the process $Z_s$ is
highly irregular (e.g., stochastically discontinuous everywhere
except at $s=1$), which is a manifestation of asymptotically fast
oscillations of the process $\Delta_n(x)$ (as well as
$\widetilde\Delta_n(u)$) in the vicinity of each point $x\in(0,2)$
(respectively, $u\in(-2,2)$).
\end{remark}

In conclusion of this section, let us comment on the asymptotics
of the random function $\Delta_n(x)$ at the ends of the limit
spectrum, that is, for $x=0$ and $x=2$. By the definition
(\ref{eq:tilde}),
$\Delta_n(0)=\lambda(0)-\sqrt{n}\,\omega(0)=\lambda_1-2\myp\sqrt{n}$,
and according to (\ref{eq:Airy})
$$
\lim_{n\to\infty}P_n\left\{\frac{\Delta_n(0)}{n^{1/6}}\le
z\right\}=F_{1}(z),
$$
where $F_1(\cdot)$ is the Tracy--Widom distribution (see Sect.~1).
However, the limit distribution of
$\Delta_n(2)=\lambda(2\sqrt{n}\mypp)$ proves to be discrete.

\begin{theorem}\label{th3}
Under the Plancherel measure\/ $P_n$, for any\/ $z\ge0$
$$
\lim_{n\to\infty} P_n\bigl\{\Delta_n(2)\le z\bigr\}=
F_{i}(0),\qquad i=[z]+1,
$$
where\/ $F_i(\cdot)$ is the distribution function of the\/
$i$-th largest point in the Airy ensemble\/ \textup{(}see
$(\ref{eq:Airy})$\textup{)}.
\end{theorem}

Indeed, using the invariance of the measure $P_n$ under the
transposition $\lambda\leftrightarrow\lambda'$ (see Sect.~1), we
have, due to (\ref{eq:Airy}),
\begin{gather*}
P_n\bigl\{\Delta_n(2)\le z\bigr\}= P_n\bigl\{\lambda'_{i}<
2\sqrt{n}\bigr\}
=P_n\biggl\{\frac{\lambda_i-2\sqrt{n}}{n^{1/6}}<0\biggr\} \to
F_i(0).
\end{gather*}

\begin{remark}\label{rm7}
\normalfont
In the ``rotated'' coordinates $u,v$, a similar result holds for both edges:
$$
\lim_{n\to\infty}
P_n\bigl\{{\textstyle\frac{1}{2}}\mypp\widetilde\Delta_n(\pm2)\le
z\bigr\}=F_{i}(0),\qquad i=[z]+1.
$$
\end{remark}

\medskip
\begin{center}
\large 3. POISSONIZATION
\end{center}

\smallskip
The proof of Theorems \ref{th1} and \ref{th2} is based on a
standard poissonization technique (see, e.g., \cite{BDJ}). Let
${\mathcal P}=\cup_{n=0}^\infty {\mathcal P}_n$ be the set of
partitions of all natural numbers (as usual, it is convenient to
include here the case $n=0$, where there is just one, ``empty''
partition $\lambda_\emptyset\vdash 0$). Set
$|\lambda|:=\sum_{\lambda_i\in\lambda} \lambda_i$ and for $t>0$
define the poissonization $P^t$ of the measure $P_n$ as follows:
\begin{equation}\label{eq:P^t}
P^t(\lambda):=e^{-t}\,t^{|\lambda|}
\left(\frac{d_\lambda}{|\lambda|!}\right)^2,\qquad
\lambda\in{\mathcal P}.
\end{equation}
Formula (\ref{eq:P^t}) defines a probability measure on the set
${\mathcal P}$, since for $\lambda\in{\mathcal P}_n$ we  have
$|\lambda|=n$ and hence
$$
\sum_{\lambda\in {\mathcal P}}
P^t(\lambda)=e^{-t}\sum_{n=0}^\infty
\frac{t^{n}}{n!}\sum_{\lambda\in{\mathcal
P}_n}\frac{d_\lambda^2}{n!}=e^{-t}\sum_{n=0}^\infty
\frac{t^{n}}{n!}=1.
$$

We first prove the ``poissonized'' versions of Theorems \ref{th1}
and \ref{th2}. Let $Y_{t}(x)$ be given by formula (\ref{eq:Y})
with $n$ replaced by~$t$.

\begin{theorem}\label{th4}
For each\/ $0<x<2$, the distribution of the random variable\/
$Y_t(x)$ with respect to the measure\/ $P^t$ converges, as\/
$t\to\infty$, to the standard normal distribution\/ ${\mathcal
N}(0,1)$.
\end{theorem}

\begin{theorem}\label{th5}
In the notations of Theorem\/ $\ref{th2}$, the random vector\/
$(Y_{t}(x_0),\dots,\allowbreak Y_{t}(x_m))$ converges in
distribution\/ \textup{(}with respect to the measure\/
$P^t$\textup{)} to a Gaussian vector\/ $(Z_{s_0},\dots,\allowbreak
Z_{s_m})$ with zero mean and the same covariance matrix~$K$.
\end{theorem}

In order to derive Theorems \ref{th1} and \ref{th2} from Theorems
\ref{th4} and \ref{th5}, respectively, one can use a standard
de-poissonization method. According to formula (\ref{eq:P^t}), the
expression for $P^t$ can be viewed as the expectation of the
random measure $P_N$, where $N$ is a Poisson random variable with
parameter~$t$:
\begin{equation}\label{eq:DeP}
P^t(A)=E\,[P_N(A)]=e^{-t}\sum_{k=0}^\infty
\frac{t^{k}}{k!}\,P_k(A).
\end{equation}
Since $N$ has mean $t$ and standard deviation $\sqrt{t}$, equation
(\ref{eq:DeP}) suggests that the asymptotics of the probability
$P_n(A)$ as $n\to\infty$ can be recovered from that of $P^t(A)$ as
$t\sim n\to\infty$. More precisely, one can prove that
$$
P_n(A) \sim P^t(A), \qquad t\sim n \to\infty,
$$
provided that variations of the probability $P_k(A)$ are small in
the zone $|k-n|\le \mbox {const}\myp\sqrt n$. In the context of
random partitions, such a de-poissonization lemma was obtained by
Johansson (see \cite{BDJ}).

\bigskip
\begin{center}
\large 4. SKETCH OF THE PROOF OF THEOREM 4
\end{center}

\smallskip
Note that, in view of (\ref{eq:omega}), the statement
of Theorem \ref{th4} is equivalent to saying that for any
$z\in\RR$
\begin{equation}\label{eq:Phi}
\lim_{t\to\infty} P^t\left\{\lambda(\sqrt{t}\myp
x)-\lceil\sqrt{t}\myp x\rceil\le 2\myp\sqrt{t}\cos \theta_x +
z\sqrt{\log t\mathstrut}\,\right\}= \Phi(2\myp\theta_x z),
\end{equation}
where $\Phi(\cdot)$ is the distribution function of the normal law
${\mathcal N}(0,1)$. Using the Frobenius coordinates
$\lambda_i-i$, set
$$
{\mathcal D}(\lambda):=\cup_{i=1}^\infty \,\{\lambda_i-i\},\qquad
\lambda\in {\mathcal P}.
$$
Consider the semi-infinite interval
$I_t:=[2\myp\sqrt{t}\cos\theta_x + z\sqrt{\log
t\mathstrut},\infty)$ and let $\# I_t$ be the number of points
$\lambda_i-i\in {\mathcal D}(\lambda)$ contained in $I_t$. Using
that the sequence $\lambda_i-i$ is strictly decreasing and
recalling the definition (\ref{eq:lambda}) of the function
$\lambda(\cdot)$, it is easy to see that relation (\ref{eq:Phi})
is reduced to
\begin{equation}\label{eq:I->Phi}
\lim_{t\to\infty}P^t\{\lambda\in {\mathcal P}:\#I_t\le \lceil
\sqrt{t}\myp x\rceil\}=\Phi(2\myp\theta_x z).
\end{equation}

For $k=1,2,\dots$,  define the $k$-point correlation functions by
$$
\rho^{\,t}_k(x_1,\dots,x_k):=P^t\bigl\{\lambda\in {\mathcal P}:
x_1,\dots,x_k\in {\mathcal D}(\lambda)\bigr\}\qquad(x_i\in\ZZ,\ \
x_i\ne x_j).
$$
The key fact is that the correlation functions $\rho^{\,t}_k$ have
a determinantal structure (see \cite{BOO,J}):
$$
\rho^{\,t}_k(x_1,\dots,x_k)=\det[J(x_i,x_j;t)]_{1\le i,j\le k}\,,
$$
with the kernel $J$ of the form
$$
 J(x,y;t)= \left\{
 \begin{array}{ll}
 \displaystyle
 \sqrt{t}\;\frac{J_xJ_{y+1}-J_{x+1}J_{y
}}{x-y}\,,& x\neq y,\\[.8pc]
 \displaystyle
 \sqrt{t}\left(J'_xJ_{x+1}-J'_{x+1}J_x\right),&x=y,
\end{array}
\right.
$$
where $J_m=J_m\bigl(2\sqrt{t}\bigr)$ is the Bessel function of
integral order~$m$.

In this situation, one can apply Soshnikov's theorem \cite{S}
(generalizing an earlier result by Costin and Lebowitz \cite{CL}),
stating that the random variable $\# I_t$
%, under some mild technical conditions,
satisfies the central limit theorem:
\begin{equation}\label{eq:CL}
  \frac{\# I_t-E[\#I_t]}{\sqrt{\mathstrut\Var[\# I_t]}}
  \stackrel{d}{\longrightarrow}{\mathcal N}(0,1)\qquad
  (t\to\infty),
\end{equation}
provided that $\Var[\#I_t]\to\infty$. Thus, in order to derive
(\ref{eq:I->Phi}) from (\ref{eq:CL}), it remains to obtain the
asymptotics of the first two moments of the random variable
$\#I_t$. The next lemma is the main technical (and most difficult)
part of the work.

\begin{lemma}
Let\/ $E^t$ and\/ $\Var^t$ denote expectation and variance,
respectively, under the measure\/ $P^t$. Then, as\/ $t\to\infty$,
\begin{align*}
E^t\bigl[\#I_t\bigr]&=\sqrt{t}\myp
x-\frac{z\myp\theta_x}{\pi}\sqrt{\log
t\mathstrut}+O(1),\\
\Var\nolimits^t\bigl[\#I_t\bigr]&=\frac{\log t}{4\pi^2}\,(1+o(1)),
\end{align*}
\end{lemma}

The proof of Lemma 1 is based on a direct asymptotic analysis of
the expressions for the expectation and variance. In so doing, the
calculations are quite laborious and heavily use the asymptotics
of the Bessel function $J_m\bigl(2\sqrt{t}\bigr)$ in various
regions of variation of the parameters.

\medskip
Finally, note that the proof of Theorem \ref{th5} follows similar
ideas using Soshnikov's central limit theorem for linear
statistics of the form $\sum_i\alpha_i\#I_{t_i}$ \cite{S2}.

\bigskip
\begin{center}
\large 5. LINK WITH KEROV'S RESULT
\end{center}

\smallskip
Let us comment on the link between our results and the limit
theorem by Kerov \cite{K} (see Sect.~1). In particular, our goal
is to explain heuristically the mechanism of the effects that take
place for the process $\widetilde\Delta_n$.

Note that if $|u-u'|= n^{-s/2}$, $0\le s\le 1$, then
$$
  s=\frac{-2\log |u-u'|}{\log n}\,.
$$
That is to say, ``time'' $s$ indexing the components of the limit
Gaussian vector in Theorem \ref{th2'}$'$, has the meaning of the
logarithmic distance between the points $u$ and $u'$, normalized
by $\log n$. From this point of view, Theorem \ref{th2'}$'$
implies that
$$
\Cov\bigl(\widetilde Y_n(u),\widetilde Y_n(u')\bigr)\approx
s_n(u,u'):=\min\left\{\frac{-2\log |u-u'|}{\log n}\,,1\right\}.
$$
In fact, in the course of the proof of Theorems \ref{th2} and
\ref{th2'}$'$ we obtain that for any $\varepsilon>0$ there exist
constants $C_1,\myp{}C_2>0$ such that for sufficiently large $n$,
the following estimate holds uniformly in $u,u'\in(-2,2)$:
\begin{equation}\label{eq:Kbound}
\Cov\bigl(\widetilde\Delta_n(u),\widetilde\Delta_n(u')\bigr)\le
\left\{\begin{array}{ll} \displaystyle
-C_1\log |u-u'|,& |u-u'|\le\varepsilon,\\
\displaystyle \hphantom{-}C_2,& |u-u'|\ge\varepsilon.
\end{array}\right.
\end{equation}

Consider now the integral of $\widetilde\Delta_n(u)$ with respect
to a test function $\varphi$:
$$
\widetilde\Delta_n[\varphi]:= \int_{-2}^2
\widetilde\Delta_n(u)\,\varphi(u)\,du.
$$
Using (\ref{eq:Kbound}) we have
\begin{align*}
\Var\bigl(\widetilde\Delta_n[\varphi]\bigr) &=
\int_{-2}^2\int_{-2}^2 \varphi(u)
\,\varphi(u')\Cov(\widetilde\Delta_n(u),\widetilde\Delta_n(u'))\,du\,du'\\[.3pc]
&\le - C_1 \iint_{|u-u'|\le \varepsilon}
\varphi(u)\,\varphi(u')\log|u-u'|\,du\,du'
\\[.3pc]
&\quad+C_2\iint_{|u-u'|\ge \varepsilon} \varphi(u)
\,\varphi(u')\,du\,du' <\infty,
\end{align*}
since the function $\log|u-u'|$ is integrable at zero. Therefore,
$\widetilde\Delta_n[\varphi]$ is bounded in distribution as
$n\to\infty$, which helps understand why Kerov's result holds
without any normalization (see Sect.~1).

\begin{remark}\label{rm8}
\normalfont
We believe that by sharpening the asymptotic estimates (\ref{eq:Kbound}), it
may be feasible to compute the limit of the variance
$\Var(\widetilde\Delta_n[\varphi])$ and thus recover the Kerov theorem
directly from the analysis of the correlation structure. We will
address this issue elsewhere.
\end{remark}

Conversely, the limiting process $\widetilde\Delta(u)$ defined in
(\ref{eq:genDelta}) can be used to get the information contained
in Theorems \ref{th1'}$'$ and \ref{th2'}$'$ (at least
heuristically). To this end, observe that since the number of
terms in a typical partition is close to $2\sqrt{n}$ (see
Remark~\ref{rm1}), it is reasonable to think that the number of
``degrees of freedom'' of a random partition $\lambda\vdash n$ is
of order of $m\asymp \sqrt{n}$, and hence the random variable
$\widetilde\Delta_n(u)$, $u\in(-2,2)$, may be represented by the
partial sum of the series (\ref{eq:genDelta})
$$
S_m(u)|_{u=2\cos\theta}:=\frac{2}{\pi}\sum_{k=2}^{m} \frac{X_k
\sin (k\theta)}{\sqrt{k}}\,,\qquad \theta\in(0,\pi).
$$
Note that for any $u=2\cos\theta$, $u'=2\cos\theta'$
\begin{align}
\notag
\Cov\bigl(S_m(u),S_m(u')\bigr)&=\frac{4}{\pi^2}\sum_{k=2}^m
\frac{\sin(k\theta)\sin
(k\theta')}{k}\\
\label{eq:CovS} &=\frac{2}{\pi^2}\sum_{k=2}^m \frac{\cos
(k(\theta-\theta'))}{k} -\frac{2}{\pi^2}\sum_{k=2}^m \frac{\cos
(k(\theta+\theta'))}{k}\,.
\end{align}
The second sum in (\ref{eq:CovS}) converges for all $\theta,\theta'\in(0,\pi)$.
For $\theta=\theta'$, from (\ref{eq:CovS}) we get
$$
\Var[S_m(u)]\sim\frac{2}{\pi^2}\sum_{k=2}^m \frac{1}{k}\sim
\frac{2\log m}{\pi^2}\sim \frac{\log n}{\pi^2}\,,
$$
and it follows that
$$
\frac{\pi S_m(u)}{\sqrt{\log n}}\stackrel{d}{\longrightarrow}{\mathcal
N}(0,1)\qquad (n\to\infty),
$$
which is in agreement with Theorem \ref{th1'}$'$. Moreover, if
$u'-u\asymp n^{-s/2}$ (and hence $\theta'-\theta\asymp n^{-s/2}$)
then the first sum in (\ref{eq:CovS}) is approximated by the
integral
\begin{align*}
\int_{2}^m\frac{\cos (x(\theta'-\theta))}{x}\,dx&=
\int_{2|\theta'-\theta|}^{m |\theta'-\theta|}\frac{\cos y}{y}\,dy
\sim
\int_{2|\theta'-\theta|}^{\varepsilon}\frac{1}{y}\,dy\\[.3pc]
&\sim
-\log |\theta'-\theta|
\sim\frac{s}{2}\log n.
\end{align*}
Hence
$$
\Cov\bigl(S_m(u),S_m(u')\bigr)\sim \frac{s\log n}{\pi^2} \qquad
(n\to\infty),
$$
as predicted by Theorem~\ref{th2'}$'$.

\begin{remark}\label{rm9}
\normalfont As already mentioned (see Remarks \ref{rm3} and
\ref{rm5} and also a comment after formula (\ref{eq:Airy})), there
is similarity between the asymptotic properties of the spectra of
random partitions and random matrices from the GUE. Our discussion
suggests that the relationship between the generalized type
convergence \cite{J0} and the localized central limit theorem
\cite{G} in the GUE can also be explained using the correlation
structure of the spectrum. One can expect that similar
relationship may be in place for other classes of random matrices
and for more general determinantal random ensembles, but this
issue needs to be studied further.
\end{remark}

\begin{remark}\label{rm10}
\normalfont Let $\widetilde Y(u)$, $u\in[-2,2]$, be a random
process with independent values, such that $\widetilde Y(u)$ has a
standard normal distribution for each $u\in(-2,2)$, and
$\widetilde Y(\pm2)=0$ (a.s.). Our results (see Theorem
\ref{th2}$'$ and Remark \ref{rm7} after Theorem~\ref{th3}) imply
that the random process $\widetilde Y_n(\cdot)$ (see
(\ref{eq:tildeY})) converges to $\widetilde Y(\cdot)$ in the sense
of finite dimensional distributions. A natural question may arise
as to whether this can be extended to weak convergence. However,
it is easy to see that the answer is negative, at least under the
natural choice of the space of continuous functions $C[-2,2]$,
because the necessary condition of tightness breaks down
(see~\cite{B}). Indeed, for any $\delta>0$, $\varepsilon>0$ and
all $u,u'\in(-2,2)$ such that $|u-u'|\le\delta$, we have
\begin{align*}
%P_n\biggl\{\sup_{|u-u'|\le\delta}|\widetilde Y_n(u)-\widetilde
%Y_n(u')|\ge\varepsilon\biggr\}&\ge
\lim_{n\to\infty}P_n\bigl\{|\widetilde Y_n(u)-\widetilde
Y_n(u')|\ge \varepsilon\bigr\} =P\bigl\{|\widetilde
Y(u)-\widetilde Y(u')|\ge \varepsilon\bigr\}
%=2\bigl(1-\Phi(\varepsilon/\sqrt{2})\bigr)
>0.
\end{align*}
%since the random variable $Z(u)-Z(u')$
%has normal distribution ${\mathcal N}(0,2)$.
Analogous remark applies to the process $Y_n(x)$, $x\in[0,2]$,
considered in the space $D[0,2]$ of right continuous functions
%\emph{cadlag}
with left limits.
\end{remark}

\bigskip
\begin{center}
\large ACKNOWLEDGMENTS
\end{center}

\smallskip
This work was done when Z.G.~Su was visiting the University of
Leeds (United Kingdom) under a Royal Society International
Fellowship financially supported by the K.C.~Wong Educational
Foundation. His research was also partially supported by the
National Science Funds of China, Grant No.~10371109. The authors
are grateful to A.~Okounkov and Ya.G.\ Sinai for their interest in
this work and also to A.V.\ Gnedin, J.T.\ Kent and Yu.V.\
Yakubovich for the discussions and helpful remarks.

\vfill\eject

\end{document}